\documentclass[12pt]{article}
 \usepackage[latin1]{inputenc}

 \usepackage{amsmath}
 \usepackage{amsfonts}
\usepackage[frenchb]{babel}
\usepackage{float}

\begin{document}

\begin{center}
{\Large
	{\sc  Planification d'expériences séquentielle dans un contexte de méta-modélisation multi-fidélité. 
                Sequential experimental design in a multi-fidelity computer experiment framework. \\}
}
\bigskip

 Loic Le Gratiet$^{1,2}$ \& Claire Cannaméla$^{2}$ \& Josselin Garnier$^{1}$
\bigskip

{\it
$^{1}$ Université Paris Diderot- Paris VII 75205 Paris Cedex 13 ; loic.le-gratiet@cea.fr
 
$^{2}$ CEA, DAM, DIF, F-91297 Arpajon, France ; claire.cannamela@cea.fr
}
\end{center}
\bigskip

{\bf R\'esum\'e.} Les gros codes de calcul  sont souvent utilisés en ingénierie pour étudier des systèmes physiques. Cependant, les simulations peuvent parfois être coûteuses en temps de calcul. Dans ce cas, une approximation de la relation entrée/sortie du code est souvent faite à l'aide d'un méta-modèle. 
En fait, un code de calcul peut souvent être lancé à différents niveaux de complexité et une hiérarchie de niveaux de code peut donc être obtenue. Il peut s'agir par exemple d'un modèle éléments finis ayant un maillage plus ou moins fin. L'objectif de nos travaux est d'étudier l'utilisation de plusieurs niveaux de code pour prédire la sortie d'un code coûteux. Le méta-modèle multi-niveaux présenté ici est un cas particulier du co-krigeage.
Après avoir détaillé le modèle de co-krigeage utilisé et sa mise en place, nous concentrerons notre présentation sur les stratégies de planification d'expériences séquentielles. En effet, un avantage du co-krigeage  est qu'il fournit, au travers de la variance de co-krigeage, une estimation de l'erreur de modèle en chaque point de l'espace des paramètres d'entrée. Ainsi, pour améliorer la précision du méta-modèle on peut enrichir séquentiellement notre base d'apprentissage aux points où la variance est la plus grande. Cependant, dans un cadre de multi-fidélité, il faut également choisir sur quel niveau de code on veut connaître la réponse. Nous présenterons ici différentes stratégies pour choisir ce niveau basées sur un résultat original qui donne la contribution de chaque code à la variance de prédiction du modèle.
\smallskip

{\bf Mots-cl\'es.} co-krigeage, méta-modélisation multi-fidélité, plan d'expériences, stratégie séquentielle.
\bigskip\bigskip

{\bf Abstract.}Large computer codes are widely used in engineering to study physical systems. Nevertheless, simulations can sometimes be time-consuming. In this case, an approximation of the code input/output relation is made using a metamodel.
Actually, a computer code can often be run at different levels of complexity and a hierarchy of levels of code can hence be obtained. For example, it can be a finite element model with a more or less fine mesh. The aim of our research is to study the use of several levels of a code to predict the output of a costly computer code. The presented multi-stage metamodel is a particular case of co-kriging which is a well-known geostatistical method. 
We first describe the construction of the co-kriging model and we focus then on a sequential experimental design strategy. Indeed, one of the strengths of  co-kriging is that it provides through the predictive co-kriging variance an estimation of the model error at each point of the input space parameter. Therefore, to improve the surrogate model we can sequentially add  points in the training set at locations where the predictive variances are the largest ones. Nonetheless, in a multi-fidelity framework, we also have to choose which level of code we have to run. We present here different strategies to choose this level. They are based on an original result which gives the contribution of each code on the co-kriging variance.
\smallskip


\section{Introduction}

Le krigeage est une classe particulière de méta-modèle qui fait l'hypothèse \emph{a priori} que la sortie que l'on essaye d'approcher est une réalisation d'un processus Gaussien. On se concentre ici  sur ce modèle et son extension pour les sorties à réponses multiples. Le lecteur pourra se référer aux livres de Santner, Williams et Notz (2003) et Rasmussen et Williams (2006) pour plus de détails sur le krigeage.
Dans le cadre des simulateurs multi-fidélités, nous possédons un code qui peut être lancé à plusieurs niveaux de précision. Ces niveaux pouvant par exemple correspondre à différentes tailles de mailles pour un code éléments finis. L'objectif de la méta-modélisation sera de construire une surface de réponse du code en utilisant une ou plusieurs de ses versions dégradées.\\

Un premier méta-modèle multi-niveaux a été mis en place par Kennedy et O'Hagan (2000) en utilisant une relation autorégressive d'ordre 1 entre deux niveaux successifs. Ce modèle est un cas particulier du co-krigeage qui est une méthode communément utilisée en géostatistique. Ensuite, Qian et Wu (2008) ont proposé une approche Bayésienne du modèle suggéré par Kennedy et O'Hagan (2000) et Forrester, Sobester et Keane (2007) ont présenté l'utilisation du co-krigeage dans un cadre d'optimisation (basée sur l'algorithme EGO : Efficient Global Optimization). Enfin, certains points limitant de la méthode ont été résolus dans le papier Le Gratiet (2011) qui présente entre autres une procédure d'estimation efficace des paramètres, une réduction de la complexité du modèle et une nouvelle approche Bayésienne évitant une implémentation prohibitive.  Dans les papiers ci-dessus, différentes méthodes de planification d'expériences ont été proposées mais aucun auteur ne s'est intéressé à la planification d'expériences séquentielle. Pourtant la hiérarchie de codes disponible peut nous permettre de mettre en place des stratégies intéressantes de planification séquentielle. La question sera la suivante : si nous voulons lancer une simulation en un nouveau point, sur quel niveau de code allons-nous la lancer ? Il va donc falloir trouver le meilleur compromis coût / précision et nous allons pour cela utiliser un résultat utile déduit de  Le Gratiet (2011).

\section{Construction du méta-modèle multi-fidélité.}

Nous présentons ici une nouvelle modélisation multi-niveaux qui est équivalente à celle proposée par  Kennedy et O'Hagan (2000)  et qui se construit de manière récursive.

\subsection{Rappels des équations du modèle AR(1) étendues au cas $s$ niveaux.}

Supposons que nous ayons $s$ niveaux de code modélisés par des processus Gaussiens $(Z_t(x))_{t=1,\dots,s}$, $x \in Q$, classés par ordre croissant de précision et tels que :

\begin{equation}\label{eq59}
\left\{
\begin{array}{l}
Z_t(x) = \rho_{t-1}(x)Z_{t-1}(x)+\delta_t(x) \\
 Z_{t-1}(x) \perp \delta_t(x) \\
\rho_{t-1}(x) =g_{t-1}^T(x) \beta_{\rho_{t-1}} \\
\end{array}
\right.
\end{equation}

où l'on définit $\forall t = 2,\dots,s$ :
\begin{equation}\label{eq60}
\delta_t(x) \sim \mathcal{PG} (f_t^T(x) \beta_t, \sigma_t^2 r_t(x,x'))
\end{equation}
et : 
\begin{equation}\label{eq61}
Z_1(x) \sim \mathcal{PG} (f_1^T(x) \beta_1, \sigma_1^2 r_1(x,x'))
\end{equation}

On fait donc ici l'hypothèse \emph{a priori} que les sorties des codes sont  des réalisations de processus Gaussiens. Le méta-modèle sera construit en conditionnant ces réalisations par les sorties connues des codes.\\

La modèlisation précédente issue de l'article de Kennedy et O'Hagan (2000) définit la relation entre les niveaux de codes. Elle est déduite de l'hypothèse suivante: $\forall x$,  si on se donne la valeur de $Z_{t-1}(x)$, on ne peut rien apprendre de plus sur $Z_t(x)$ à partir des autres résultats $Z_{t-1}(x')$ pour $x \neq x'$.\\

Soit $\mathcal{Z} = (\mathcal{Z}_1^T,\dots,\mathcal{Z}_s^T)^T$ le vecteur Gaussien contenant le valeurs des processus $(Z_t(x))_{t=1,\dots,s}$ aux points des plans d'expériences $(D_t)_{t=1,\dots,s}$ avec  $D_s \subseteq D_{s-1} \subseteq \dots \subseteq D_1$. Si on note $\beta = (\beta_1^T,\dots,\beta_s^T)^T$, $\beta_\rho = (\beta_{\rho_1}^T,\dots,\beta_{\rho_{s-1}}^T)^T$, $\sigma^2 = (\sigma_1^2,\dots,\sigma_s^2)$  et $z  = (z^1,\dots,z^s)$ les réponses connues des codes, on a :

\begin{displaymath}
\forall x \in Q \qquad [Z_s(x)|\mathcal{Z}=z,\beta,\beta_\rho,\sigma^2] \sim \mathcal{N}\left(m_{Z_s}(x),s_{Z_s}^2(x)\right)
\end{displaymath}

où :
\begin{equation}\label{eq62}
m_{Z_s}(x) = h_s'(x)^T\beta + t_s(x)^TV_s^{-1}(z-H_s\beta)
\end{equation}

et :
\begin{equation}\label{eq63}
s_{Z_s}^2(x) = \sigma_{Z_s}^2(x) - t_s(x)^T V_s^{-1} t_s(x)
\end{equation}

La moyenne de co-krigeage $m_{Z_s}(x)$ constituera le méta-modèle  sur le niveau le plus précis  construit à partir des réponses connues sur les $s$ niveaux de code et la variance de co-krigeage $s_{Z_s}^2(x)$ donnera une estimation de l'erreur de modèle.\\

La matrice $V_s$ représente la matrice de covariance du vecteur Gaussien  $\mathcal{Z}$, le vecteur $t_s(x)$ est le vecteur de covariance entre $Z_s(x)$ et  $\mathcal{Z}$,  $H_s \beta$ est la moyenne de  $\mathcal{Z}$, $h_s'(x)^T\beta$ est la moyenne de $Z_s(x)$ et  $\sigma_{Z_s}^2(x)$ sa variance. Tous ces éléments sont construits à partir du vecteur d'expériences au niveau $t$  et de la covariance entre $Z_t(x)$ et $Z_{t'}(x')$ explicités ci-dessous :
\begin{equation}\label{eq69}
h_t'(x)^T = \left( \left(\prod_{i=1}^{t-1} {\rho_{i}(x) }\right)f_1^T(x) , \left(\prod_{i=2}^{t-1} {\rho_{i}(x) }\right)f_2^T(x) , \dots , \rho_{t-1}(x) f_{t-1}^T(x) , f_t^T(x) \right)
\end{equation}

et pour $t > t'$ :
\begin{equation}\label{eq70}
\mathrm{cov}(Z_t(x),Z_ {t'}(x')|\sigma^2,\beta,\beta_\rho) =
\left(\prod_{i=t'}^{t-1}\rho_i(x)\right) \mathrm{cov}(Z_{t'}(x),Z_{t'}(x')|\sigma^2,\beta,\beta_\rho)
\end{equation}
avec :
\begin{equation}\label{eq70}
\mathrm{cov}(Z_t(x),Z_t(x')|\sigma^2,\beta,\beta_\rho) =
\sum_{j=1}^{t}{\sigma_{j}^2\left( \prod_{i=j}^{t-1} {\rho_{i}^2(x) } \right)r_j(x,x')}
\end{equation}

\subsection{Ecriture récursive du modèle AR(1)}

Nous allons ici présenter une nouvelle modélisation multi-niveaux qui est en fait équivalente à la précédente. Considérons la modélisation suivante pour $t=2,\dots,s$ :

\begin{equation}\label{eq71}
\left\{
\begin{array}{l}
Z_t(x) = \rho_{t-1}(x)[Z_{t-1}(x)|\mathcal{Z}_{t-1}=z^{t-1}]+\delta_t(x) \\
 Z_{t-1}(x) \perp \delta_t(x) \\
\rho_{t-1}(x) = g_{t-1}^T(x) \beta_{\rho_{t-1}} \\
\end{array}
\right.
\end{equation}

avec  $D_s \subseteq D_{s-1} \subseteq \dots \subseteq D_1$. La seule différence avec la précédente modélisation est que l'on exprime $Z_t(x)$ (le processus Gaussien modélisant le code de niveau $t$) en fonction du processus $Z_{t-1}(x)$ conditionné par ses réponses connues $z^{t-1}$ sur le plan $D_{t-1}$. Les $(\delta_t(x))_{t=2,\dots,s}$ sont définis comme précédemment et on a  pour tout $t=2,\dots,s$ :
\begin{equation}\label{eq75}
\left[Z_{t}(x)|\mathcal{Z}_{t}=z^{t}\right] \sim
\mathcal{N}
\left(
\mu_{Z_{t}}(x) ,
s_{Z_{t}}^2(x) 
\right)
\end{equation}

où :
\begin{equation}\label{eq76}
\mu_{Z_{t}}(x) = \rho_{t-1}(x) \mu_{Z_{t-1}}(x)  + f_t^T(x)\beta_t + 
r_t^T(x) R_t^{-1}\left(
z^t - \rho_{t-1}(D_t)\odot z_{t-1}(D_t) - F_t \beta_t
\right)
\end{equation}

et:
\begin{equation}\label{eq77}
s_{Z_{t}}^2(x)  =  \rho_{t-1}^2(x)s_{Z_{t-1}}^2(x) + \sigma_t^2\left(
1-r_t^T(x) R_t^{-1}r_t(x)
\right)
\end{equation}

Le symbole $\odot $ représente le produit matriciel élément par élément. $R_t$ est la matrice de correlation de  $D_t$  avec le noyau $r_t(x,x')$ et $r_t^T(x)$ est le vecteur de covariance entre $x$ et $D_t$ avec ce même noyau. On note $\rho_t(D_{t-1})$ le vecteur contenant les valeurs de  $\rho_t(x)$ pour  $x\in D_{t-1}$, $z_t(D_{t-1})$ celui contenant les sorties connues de $Z_t(x)$ sur $D_{t-1}$  et  $F_t$ est la matrice d'expériences contenant les valeurs de $f_t(x)^T$ sur $D_t$. \\

La moyenne $\mu_{Z_{t}}(x)$ constitue le méta-modèle sur le code de niveau $t$ sachant les réponses connues des $t$ premiers niveaux de code et la variance $s_{Z_{t}}^2(x)$ représente l'erreur de ce modèle. La moyenne et la variance de co-krigeage sur le niveau $t$ s'exprimant à partir de celles du niveau $t-1$, on a bien une expression récursive du modèle proposé par \cite{KO00}. La modélisation proposée ici est puissante car elle montre que construire un co-krigeage à $s$ niveaux est équivalent à construire $s$ krigeages indépendants.

\section{Planification d'expériences séquentielle pour le co-krigeage}

La modélisation récursive présentée précédemment va nous permettre de construire une méthode simple de planification d'expériences séquentielle. Nous travaillerons sur un des principaux critères de planification séquentielle utilisés dans le krigeage : le MSE (Mean Squarred Error). La stratégie consiste simplement à trouver le point maximisant la variance de prédiction s'écrivant :

\begin{displaymath}
s^2_{Z_t}(x)= \rho_{t-1}^2(x)s^2_{Z_{t-1}}(x)+\sigma_t^2\left(1-r_t^T(x)R_t^{-1}r_t(x)\right)  \qquad  \forall t=2,\dots,s
\end{displaymath}
\begin{displaymath}
s^2_{Z_1}(x)=\sigma_1^2\left(1-r_1^T(x)R_1^{-1}r_1(x)\right) 
\end{displaymath}

Soit $\tilde{x}=  \mathrm{ argmax}_x \, s^2_{Z_s}(x)$. On a alors $\forall t=1,\dots,s$ :

\begin{displaymath}
s^2_{Z_t}(\tilde{x})= \rho_{t-1}^2(\tilde{x})s^2_{Z_{t-1}}(\tilde{x})+\sigma_t^2\left(1-r_t^T(\tilde{x})R_t^{-1}r_t(\tilde{x})\right) 
\end{displaymath}

Supposons que l'on décide de lancer le code $t-1$ en $\tilde{x}$, on a alors $s^2_{Z_{t-1}}(\tilde{x})=0$ et :

\begin{displaymath}
s^2_{Z_t}(\tilde{x})= \sigma_t^2\left(1-r_t^T(\tilde{x})R_t^{-1}r_t(\tilde{x})\right) 
\end{displaymath}

$\rho_{t-1}^2(\tilde{x})s^2_{Z_{t-1}}(\tilde{x})$ représente donc la part de la variance due au code de niveau $t-1$ et $\sigma_t^2\left(1-r_t^T(\tilde{x})R_t^{-1}r_t(\tilde{x})\right) $ représente la part  de la variance due au code de niveau $t$. Une stratégie naturelle de planification d'expériences est alors de trouver le point $\tilde{x}$ maximisant la variance, de regarder la part de chaque niveau de code sur $s^2_{Z_s}(\tilde{x})$ et de lancer une simulation sur le niveau de code le plus intéressant (selon un critère coût/précision qui dépendra du cas d'application). \\

Il est donc nécessaire de  construire un critère qui nous permette de déterminer quel niveau de code est le plus intéressant. Nous illustrons dans ce résumé, une méthode pour choisir ce niveau dans le cas $s=2$. On introduit pour cela la moyenne de la variance de prédiction (où $\mu(x)$ représente la mesure de la loi des entrées) :

\begin{displaymath}
\mathrm{imse} = \int_Q s^2_{Z_2}(x) \, d \mu(x)
\end{displaymath}

Considérons $\tilde{x}=  \mathrm{ argmax}_x \, s^2_{Z_2}(x)$, dans le cas à 2 niveaux nous pouvons  lancer soit $z_1(\tilde{x})$ soit $z_2(\tilde{x})$. Supposons que l'on connaisse $z_1(\tilde{x})$, on aura alors :

\begin{displaymath}
s^2_{Z_2}(\tilde{x}) = \sigma_2^2\left(1-r_2^T(\tilde{x})R_2^{-1}r_2(\tilde{x})\right) 
\end{displaymath}

Si $s^2_{Z_2}(\tilde{x}) < \mathrm{imse}$ cela signifie que l'erreur de modèle en $\tilde{x}$ est inférieure à l'erreur  de modèle moyenne.  Il n'est donc pas nécessaire de lancer le code de niveau 2 en ce point car cela n'apportera pas en moyenne  plus d'information qu'un autre point. En revanche, si $s^2_{Z_2}(\tilde{x}) > \mathrm{imse}$, cela signifie que l'erreur de modèle en $\tilde{x}$ est plus importante que l'erreur moyenne et il sera donc avantageux de lancer le code de niveau 2 en ce point.

\label{fin}

\begin{thebibliography}{20}
\bibitem[1]{AF07} \begin{sc} Forrester, A. I. J., Sobester, A.  \& Keane,  A. J.\end{sc} (2007), \textit{Multi-fidelity optimization via surrogate modelling},  Proc. R. Soc. A \textbf{463}, 3251-3269.\\

\bibitem[2]{KO00} \begin{sc}  Kennedy, M. C. \& O'Hagan, A.\end{sc} (2000), \textit{Predicting the output from a complex computer code when fast approximations are available}, Biometrika  \textbf{87}, 1-13.\\


\bibitem[3]{LLG11} \begin{sc} Le Gratiet, L.  \end{sc}(2011), \textit{Bayesian analysis of hierarchical multi-fidelity codes}, Arxiv preprint arXiv:1112.5389, 2011.\\

\bibitem[4]{R06} \begin{sc} Rasmussen, C. E.  \&  Williams, C. K. I.\end{sc} (2006), \textit{Gaussian Processes for Machine Learning}, the MIT Press.\\

\bibitem[5]{S03} \begin{sc} Santner, T. J., Williams, B. J.  \&  Notz, W. I. \end{sc} (2003),  \textit{The Design and Analysis of Computer Experiments}, New York: Springer.\\

\bibitem[6]{QW07} \begin{sc} Qian, Z. \& Jeff Wu, C. F.\end{sc} (2008),  \textit{Bayesian Hierarchical Modeling for Integrating Low-accuracy and High-accuracy Experiments}, Technometrics \textbf{50}, 192-204.\\
\end{thebibliography}
\end{document}